\documentclass{article}
\usepackage{amsmath,amssymb,amsfonts}

\newtheorem{conjecture}{Conjecture}
\newtheorem{proposition}{Proposition}

\begin{document}

\title{On the volume of the polytope of doubly stochastic matrices}
\author{Clara S. Chan \thanks{clara@ccr-p.ida.org} \and 
        David P. Robbins \thanks{robbins@ccr-p.ida.org}}
\maketitle

\begin{abstract}
We study the calculation of the volume of the polytope $B_n$ of
$n \times n$ doubly stochastic matrices; that is, the set of real non-negative
matrices with all row and column sums equal to one.

We describe two methods.  The first involves a decomposition of the
polytope into simplices.  The second involves the enumeration of
``magic squares'', {\it i.e.,\/}\
$n \times n$ non-negative integer matrices whose
rows and columns all sum to the same integer.

We have used the first method to confirm the previously known values through
$n=7$. This method can also be used to compute the volumes of faces
of $B_n$.  For example we have observed that the volume of
a particular face of $B_n$ appears to be a product of Catalan numbers.

We have used the second method to find the volume for $n=8$, which we believe
was not previously known.
\end{abstract}

\section{Introduction}
\label{sec:intro}
We study the calculation of the volume of the polytope $B_n$ of
$n \times n$ doubly stochastic matrices; that is, the set of real nonnegative
matrices with all row and column sums equal to one.   This polytope is
sometimes known as the Birkhoff polytope or the assignment polytope.
We will describe and evaluate two methods for computing the volume of $B_n$.

In the first method we decompose $B_n$ into a disjoint union of simplices all
of the same volume and count the simplices.  The fact that this can be 
done appears in \cite{St1}.  This method applies to any face of $B_n$ as well.

In the second method we count the number of $n \times n$ nonnegative integer
matrices with all row and column sums equal to $t$ (sometimes called magic
squares) for suitable values of $t$.
These numbers allow us to compute the Ehrhart polynomial of $B_n$,
which (essentially) has the volume of $B_n$ as its leading coefficient.
It appears that this has been the most common method of computing the
volume of $B_n$.  Sturmfels reports in \cite{Stu} on other work in 
which the volume of $B_n$ has been computed for $n$ up to 7.
We have also used this method to compute the volume when $n=8$.

This study is largely expository since the two methods are not
new.  However, the details about how we carry out these methods
may be of interest.  We are not aware of any reports of others who
have carried out the simplicial decomposition method.  

As a byproduct of our program for carrying out the simplicial
decomposition method, we are easily able to compute the volume of any
face (of any dimension) of $B_n$ provided that $n$ is not too large.
This allowed us to
discover that a certain special face of $B_n$ has a volume which appears
to be given by a simple product formula.
This formula is given in Conjecture \ref{conj:1}.

Our study resulted from a question \cite{Mi} of Victor Miller, who asked how
one could generate a doubly stochastic matrix uniformly at random.  
It is not hard to see that it would be easy to generate
a random doubly stochastic matrix if one could easily calculate the
volume of any face of $B_n$.  However the method described here for
calculating face volumes is practical only for small $n$.

In what follows we will make use of some well known properties of the 
face structure of $B_n$: 
the vertices of $B_n$ are precisely the $n!$ $n \times n$
permutation matrices; on the other hand,
for each pair $(i,j)$ with
$1 \le i,j \le n$, the doubly stochastic matrices with
$(i,j)$ entry equal to 0 form a {\it facet\/}\ (maximal proper face) of $B_n$ and all 
facets arise in this way.  See \cite{BS} for further properties and 
references.  

In general, it is convenient to identify the faces of $B_n$ with certain
$n \times n$ matrices of 0's and 1's, as follows.

First we identify a 0-1 matrix with the set of entries in the matrix that
are 1's. Thus, for two 0-1 matrices $A$ and $B$ of the same size,
we can define their union $A \cup B$ as the 0-1 matrix whose set
of 1's is the union of the sets of 1's of $A$ and $B$. {\it e.g.,\/}\
\[
\left(\begin{array}{ccc} 1 & 0 & 0 \\ 0 & 1 & 0 \\ 0 & 0 & 1 \end{array}\right)
 \; \cup \;
\left(\begin{array}{ccc} 0 & 1 & 0 \\ 1 & 0 & 0 \\ 0 & 0 & 1 \end{array}\right)
 \; = \;
\left(\begin{array}{ccc} 1 & 1 & 0 \\ 1 & 1 & 0 \\ 0 & 0 & 1 \end{array}\right).
\]
Similarly we can speak of one 0-1 matrix containing another and so forth.

Now to each face $F$ of $B_n$, we associate the matrix
$M$ which is the union of the vertices (permutation matrices) in $F$.
The facets of $B_n$ containing $F$ are precisely those associated with the
zero entries of $M$.  Since any face is the intersection of the facets
containing it, any permutation matrix contained in $M$ must be a vertex
of $F$.
Thus the vertices of $F$ are precisely the permutation matrices contained
in $M$, so we can recover $F$ from $M$.
In this way we identify the faces of $B_n$ with the set of 0-1
matrices which are unions of permutation matrices. 
Note that not every 0-1 matrix corresponds to a face of $B_n$. For example
\[
\left(\begin{array}{cc} 0 & 1 \\ 1 & 1 \end{array}\right)
\]
is not a union of permutation matrices, hence not a face of $B_2$.

\section{Volume}
\label{sec:2}
It is easy to see that the dimension of $B_n$ is $(n-1)^2$.
Strictly speaking, the volume we wish to compute is the
$(n-1)^2$-volume of $B_n$ regarded as a subset of
$n^2$-dimensional Euclidean space.  Thus, for example, the polytope
$B_2$ consists of the line segment joining the matrices 
\[
\left(\begin{array}{cc} 1 & 0  \\ 0 & 1  \end{array}\right) \quad\hbox{and} \quad
\left(\begin{array}{cc} 0 & 1  \\ 1 & 0  \end{array}\right)
\]
and hence its volume is 2.

An $n \times n$ doubly stochastic matrix is determined
by its upper left $(n-1)\times (n-1)$ submatrix.  The set of $(n-1)\times (n-1)$
matrices obtained this way is the set $A_n$ of all nonnegative $(n-1)\times (n-1)$
matrices with row and column sums $\le 1$ such that the sum of all the
entries is at least $n-2$.  This is affinely isomorphic to
$B_n$. In the Appendix we show that the ratio of the volume of $B_n$
to the volume of $A_n$, regarded as a subset of Euclidean $(n-1)^2$ space,
is $n^{n-1}$.  In some ways the volume of $A_n$ is easier to
understand since its dimension is equal to the dimension of its ambient space.

James Maiorana \cite{Ma} (and probably others) noted a Monte Carlo method for approximating
the volume of $A_n$. Consider the set $C_n$ of $(n-1)\times (n-1)$ nonnegative
matrices with row sums (but not necessarily column sums) $\le 1$.  This is the
Cartesian product of $n-1$ unit simplices in Euclidean $(n-1)$-space so 
its volume is $ \frac{1}{(n-1)!^{n-1}}$.
It is easy to choose points in $C_n$ uniformly at random.  The probability
$\alpha_n$ that such a point is in $A_n$ is the ratio of the volume of 
$A_n$ to that of $C_n$.
Thus we can
run Monte-Carlo trials to estimate $\alpha_n$ and 
hence the volume of $A_n$.

For large $n$, this Monte Carlo method is impractical since $\alpha_n$ is too small.
However, it is useful for checking computations for small $n$.
A lower bound for $\alpha_n$ is given by Bona in \cite{Bo}.

There is a more natural unit for the volume of $B_n$ and
its faces.  This is based on the fact that the vertices of $B_n$
are integer matrices.  Suppose that $F$ is a $d$-dimensional
face of $B_n$.  Since its vertices have integer coordinates, the integer
points in the affine span of $F$ comprise a 
$d$-dimensional affine lattice $L$.  Given such a lattice there is
a minimum volume of any $d$-simplex with vertices in $L$.  Lattice points
$w_0,\dots,w_d$ are the vertices of one of these minimum volume simplices
if and only if every point of $L$ is uniquely expressible in the
form $ \sum_{i=0}^d  k_i w_i$, where the $k_i$'s are integers whose sum is 1.
The {\it relative volume} of a face $F$ is the volume of
$F$ expressed in units equal to the volume of a minimal simplex in $L$.
The relative volume of a face is the
same whether regarded as a face of $B_n$ or as a face of $A_n$, since
the mapping from $B_n$ to $A_n$ (by taking the upper left $(n-1) \times (n-1)$ minor) 
preserves integrality of points.

Here are the currently known relative volumes of $B_n$.
\[
\begin{array}{cr}
n & \mbox{Relative Volume of $B_n$} \\
1 & 1\\
2 & 1\\
3 & 3\\
4 & 352\\
5 & 4718075\\
6 & 14666561365176\\
7 & 17832560768358341943028\\
8 & 12816077964079346687829905128694016
\end{array}
\]

To convert relative volumes to true volumes, we need to know the
volume of a minimal simplex of $A_n$.  

But the affine span of $A_n$ is all of $(n-1)^2$-dimensional space.  Hence
the volume of a minimal simplex in $A_n$ is $\frac{1}{((n-1)^2)!}$,
and the volume of a minimal simplex in $B_n$ is $\frac{n^{n-1}}{((n-1)^2)!}$.

\section{Triangulations}
\label{sec:3}

We call the first method for computing the volume of $B_n$ the 
{\it triangulation method.\/}\  The method applies to the
calculation of the volume of any polytope.  The
essence is that we decompose the polytope into simplices and
sum the volumes of the simplices.

For $B_n$ we have used a standard method of decomposing a 
polytope $P$ into simplices.  See for example \cite{St1}.
To decompose $P$ into simplices, we
choose an arbitrary vertex $v$ and form the collection of facets of $P$ 
opposite $v$ (facets of
$P$ not containing $v$.)  We then
recursively triangulate each facet.  The triangulation of 
$P$ is then formed by adding our chosen vertex to each simplex in
the triangulation of each of the facets.

The standard triangulations of $B_n$ and its faces have an unusual property,
given in \cite{St1}, for which we provide a self-contained proof below.

\begin{proposition}
In any standard triangulation of a face $F$ of $B_n$, every
simplex has minimal volume in the affine lattice determined by $F$.
\end{proposition}

{\it Proof:\/}\
Let $F$ be a $d$-dimensional face of $B_n$,
$v_0$ any vertex in $F$, and $G$ a facet of $F$ opposite $v_0$.
Suppose that a simplex in a standard triangulation of $G$
has vertices $v_1,v_2,\dots,v_d$.
We need to prove that the set of integer points of the affine space determined by $F$ is 
the same as the set of points
$\sum_{i=0}^{d} k_i v_i $,
where the $k_i$ are integers whose sum is 1.  

Of course all the integer combinations are in the affine span.
The question is whether there are any other points.

Any integral point of the affine span can be uniquely expressed in the form
$\sum_{i=0}^{d} r_i v_i $,
where the $r_i$'s are real numbers with sum 1.

Since $v_0$ is not in the face $G$, there is a facet of
$B_n$ containing $G$ but not $v_0$.  Thus
$v_0$ must have at least one
entry equal to $1$ in the same position where all $v_i$, $i\ge 1$, have zeroes.
Thus, in the hypothetical combination above, $r_0$ must be an integer.
If we add $r_0(v_1-v_0)$ to the combination above, we obtain another integral
point in the affine span of $G$.  It follows, using induction, that $r_1+r_0$, and
$r_2,\dots,r_d$ are integers and therefore all the $r$'s are integers,
as desired.\\

Note that, as a corollary, in any standard triangulation
of a face of the $B_n$, the number of simplices in the
triangulation is equal to the relative volume.  

We also obtain an important computational principle.  Given a face $F$
of $B_n$ and a vertex $v$ of $F$, the relative volume
of $F$ is the sum of the relative volumes of facets of $F$ opposite
$v$.

\section{The Triangulation Method for $B_n$}
\label{sec:4}

We now describe the triangulation method for computing the volume of
$B_n$.  This is simply an elaboration of the principle that
the relative volume of a face is the sum of the relative volumes of the
faces opposite any vertex.

We apply this principle recursively.  To get started we use the fact
that the relative volume of any zero-dimensional face of $B_n$ is 1.

In the most naive plan we calculate the relative
volumes of all faces.  We first produce a list of all faces of each
dimension.  For dimension 0, we know all the relative volumes are 1.
Then, for each face $F$ of dimension $d$ we select a vertex and
find the opposite facets (of dimension $d-1$). Assuming recursively
that their relative volumes have already been computed, we now find the relative volume
of $F$ by summing the relative volumes of the facets.

There are two serious drawbacks to the naive plan. 

Perhaps the most pressing problem is that we need to compute the volumes of
an extremely large number of faces, since quite a few of the $2^{n^2}$
possible 0-1 matrices are actually faces of $B_n$.
Here we have recourse to a single important trick.  If we
permute the rows and columns of the matrix representing a face, we obtain
the matrix of another face with the same volume.  Also if we transpose
a matrix representing a face, we obtain another face of the same
volume.  We regard matrices which can be obtained from each other by
these operations as equivalent.  We can cut down on the cost of our
algorithm if we compute the volume only for a single
``canonical'' face in each equivalence class.

The next most difficult problem is to produce the lists of faces.
The most practical method that we found for producing faces is to
start with the single $(n-1)^2$-dimensional face, $B_n$ itself,
and successively produce faces of lower
dimension by intersecting with a facet of $B_n$.
While producing the faces we save the subface information so that
we can look up the volumes when we are done.  Unfortunately
we need to construct a very large partially ordered set of faces before
we can calculate any volumes since the only volumes we know
are those of the zero-dimensional faces.
While the cost in memory is not so bad for $n$ less than 8, when we reach
$n=8$, we seem to need about 200 gigabytes of intermediate storage.
If the memory were available, the computation of the volumes
would be relatively easy.  In fact we are able to carry out a
substantial fraction of the work before running out of memory.

There are two phases to our algorithm.
In the first phase we construct a collection of faces
together with information about which ones are facets of
which others.  In particular, we successively compute, for
$d=(n-1)^2,(n-1)^2-1,\dots,0$, a collection ${\cal F}_d$ of
$d$-dimensional faces of $B_n$.  We begin by
setting ${\cal F}_{(n-1)^2} = \{ B_n \}$,
{\it i.e.,\/}\ consisting of just the all 1's matrix representing
$B_n$ itself.

Given ${\cal F}_d$ we produce ${\cal F}_{d-1}$ as follows.
Start with ${\cal F}_{d-1}=\emptyset$.
For each face $F\in {\cal F}_d$ we select a vertex $v\in F$.
We then find
the facets of $F$ opposite $v$, canonicalize these faces, and add them
to ${\cal F}_{d-1}$.  Having done this for all $F\in {\cal F}_d$,
we sort ${\cal F}_{d-1}$ and remove the duplicates.
Then, for each face $f\in{\cal F}_{d-1}$, we save a list of pointers to
the faces in ${\cal F}_d$ from which $f$ arose.  
(Equivalent faces can appear several times as opposite faces
of the same face.  When this happens, we include the pointer in the list of
pointers multiple times.)

This completes the first phase.  In the
second phase we start with ${\cal F}_0$ and work up to higher
dimensions, calculating the relative volume of every saved face until
we obtain the volume of $B_n$ itself.  This is quite fast, requiring just
one addition for each saved pointer.

Note that once the pointers are constructed we do not 
need the faces themselves, unless we want to know which
face has each of the intermediate volumes we are computing.

For larger values of $n$, the accumulators used for calculating
the volumes will overflow.   But we can get around this
problem by using multiple precision arithmetic or by
performing the volume calculation several times modulo various
primes and combining the results with the Chinese Remainder Theorem.

The main computational work of our algorithm takes place in three steps.

\begin{enumerate}
\item for each face $F\in {\cal F}_d$, find a vertex $v\in F$.  
\item determine the facets of $F$ opposite $v$.
\item put these opposite facets into canonical form.
\end{enumerate}

We now describe how each of these steps is done.

One important decision is the data structure for storing faces. We
identify each face with the 0-1 matrix which is the union of its
vertices (regarded as permutation matrices).  For $n\le 8$ it is
convenient to represent each face as $n^2$ bits of a single word,
where the words of a (64-bit) computer are regarded as 64-long arrays
of bits.

In Step 1 we are given a face $f$ represented by a 0-1 matrix and
we are looking for a permutation matrix $\pi$ contained in $f$.  This
could be done with the assignment algorithm or one of the methods for
finding maximum matchings, but for the small values of $n$ that we
were using, it was quicker to use a backtracking search method, as follows.
The matrix $f$ has at least one 1 in its first row.
We guess one of these as the location of the 1 in the first row of $\pi$.
We then guess the location of the 1 in the second row of $\pi$,
bearing in mind that it cannot be in the same column as the 1 in the first row.
We continue this way searching for the location of the 1 in subsequent rows.
We backtrack if we reach a row in which there are no feasible choices.

Now we consider Step 2.  Given a face $f$ and a vertex $\pi$ we need to find
the facets of $f$ opposite $\pi$.

For a moment let us ignore $\pi$ and consider the general problem of
constructing the facets of $f$.  The main principle is that each facet
of $f$ can be obtained by intersecting $f$ with a facet of $B_n$ that does not
contain $f$. Consider the facet corresponding to the pair $(i,j)$.
The facet does not contain $f$ if $f_{ij}=1$.  To intersect $f$ with this
facet we start by replacing $f_{ij}$ with 0, obtaining a 0-1
matrix $g$.  The face which is the intersection of $f$ with the facet
$(i,j)$ is then the largest face $h$ contained in $g$.  
The matrix $h$, which is the union of the permutation matrices contained in
$g$, can be strictly contained in $g$.  Given one of the 1's in $g$, to test
whether it is in $h$, we search for a permutation matrix in $g$ which uses
the 1 in question.  This can be done with our
backtracking search algorithm.  The 1 in question is in $h$ precisely
when this search succeeds.  When the position of a 1 in $g$ is zero in
$h$ we say it is {\it forced} to zero.

For example if $n=3$ and $f$ is the 3-dimensional face with matrix
\[
\left(\begin{array}{ccc} 0 & 1 & 1 \\ 1 & 1 & 1 \\ 1 & 1 & 1 \end{array}\right)
\]
then the intersection of $f$ with the facet corresponding to the middle entry
of the top row is one-dimensional, with matrix
\[
\left(\begin{array}{ccc} 0 & 0 & 1 \\ 1 & 1 & 0 \\ 1 & 1 & 0 \end{array}\right).
\]

In this example two zeroes in the last column are forced.
This example also shows that although every facet of $f$ is the intersection of
$f$ with a facet of $B_n$, the converse is not true and the dimension of the
intersection can be too small to be a facet of $f$.

There are some additional simplifications when we search for the facets of $f$
opposite a given vertex $\pi$ of $f$.  If $g$ is a facet of $f$ not containing
$\pi$, then $g$ must contain a 0 in place of one of the 1's of $\pi$.  Thus
there are at most $n$ facets of $f$ opposite $\pi$.
Observe that if $g$ is a facet of $f$ not containing
$\pi$, then $g \cup \pi$ is a union of permutation matrices and
therefore a face of $B_n$ containing $g$ and $\pi$.  Thus $g \cup \pi=f$.
This implies an important and helpful principle.  When we introduce a 0 at a
1 of $\pi$ and this results in a facet of $f$, 
then the only other positions that might be forced to zero are
those of the other 1's of $\pi$. Thus we can loop through the $n$ 1's of $\pi$
one at a time and, for each of these, introduce a 0 and determine what other
1's of $\pi$ are forced to be 0 and produce accordingly a matrix, which we
call a candidate.  We obtain a set of
$n$ candidates among which all the facets opposite $\pi$ must occur.
(This list can have duplicates which we remove.)  Of these
candidates the facets are those which are maximal under inclusion.
Indeed, it is clear that every candidate contains a face that has
the same intersection with $\pi$.  But this face is contained in a
facet which has an intersection with $\pi$ that is at least as large.
Thus every candidate is contained in at least one facet,
and the facets are precisely the maximal candidates.

Finally we describe Step 3, which we call ``canonicalization''.

The most straightforward way to choose a canonical form 
for a face $f$ is to apply every element of our group of
symmetries to $f$ and choose the image of $f$ with the least value (where
the bit pattern $f$ is regarded as an integer.)
But this is prohibitively slow.  

Instead we make use of certain special functions, which we call scores,
which assign integers to every row and column of a 0-1 matrix.
The scores have the special property that
when rows are permuted,
the row scores are permuted the same way leaving the column scores 
unchanged, whereas, when columns are permuted, the column scores are permuted
the same way leaving the row scores unchanged.  An example of an allowable 
score is to assign 
to each row its row sum and to each column its column sum.

Given such scores we say that a matrix is in standard form if it satisfies the
following three properties:

\begin{enumerate}
\item the column scores are weakly increasing.
\item the row scores are weakly increasing.
\item in the case of tied row scores the rows are ordered lexicographically
as bit strings.
\end{enumerate}

For a given 0-1 matrix, once
its row and column scores have been computed it is easy to put a matrix
and its transpose into standard form by forcing each of the three
conditions above in the listed order.  

For each face constructed,
we put both the face and its transpose into standard form
and finally choose the smaller of these two, regarded as integers,
as the ``canonical'' form that is saved.

Note that we are abusing terminology a little here since although
the method always replaces a face by an equivalent face, it is
conceivable that equivalent faces will canonicalize to distinct faces.
When this happens, we still obtain correct volumes, but we end up doing
work which could be avoided if the equivalence were recognized.
However, if this event is rare, we obtain almost all the savings
of true canonicalization as described above (but without the excessive cost).

It turns out that just using row and column sums as the score functions
fails to recognize a substantial number of equivalences.  What we need
are scores that tend to assign different values to different rows and
columns.  Slightly more complicated scores do better.  Given a
column score, we can produce a more complicated row score by assigning
to each row the sum (or any symmetric function) of the values of the
column scores of those columns for which 1's occur in the given row.
Similarly a row score can be used to produce a more complicated column
score.  We can also add two row scores to obtain another row score or
two column scores to obtain another column score.  By combining steps like this
we produced scores that were better at
distinguishing rows (and columns) without being much more expensive to
compute.

This concludes our description of the triangulation method.
As mentioned earlier it is reasonably practical for $n < 8$.
The times required on a 500mhz DEC alpha were as follows:

\[
\begin{array}{cr}
 & \mbox{Time in seconds} \\
n<6 & \mbox{less than 0.1}\\
n=6 & 0.63\\
n=7 & 250.1
\end{array}
\]

Although the volumes of $B_n$ do not seem to follow a recognizable
pattern, it seemed conceivable that there would be faces of $B_n$ for
which the relative volumes had interesting properties.
One fairly natural class is the set of matrices for which the set of
of zeroes of the matrix form a Young tableau in a corner 
of the matrix.  

Since our triangulation method applies to any face of $B_n$, we were able
to check some natural classes of faces.  It turned out that for the 
simplest non-trivial Young tableau faces the volumes apparently
obey a simple rule, although we have not been able to supply a proof.
More precisely, suppose that $n\ge 2$ and that $F_n$ is the $n\times n$
matrix whose $(i,j)$ entry is 1 when $j \le i+1$ and 0 otherwise.
Then $F_n$ is a union of permutation matrices corresponding to a face of $B_n$
of dimension $\binom{n}{2}$ with $2^{n-1}$ vertices and we have the following\\

\begin{conjecture}\label{conj:1}
The relative volume of $F_n$ is the product
\[
\prod_{i=0}^{n-2} \frac{1}{i+1}\binom{2i}{i}
\]
of the first $n-1$ Catalan numbers.\\
\end{conjecture}

We have verified this for $n \leq 12$.\\

Finally, we give some miscellaneous observations which may be useful
but do not actually enter our algorithm.

\begin{itemize}

\item In our method, we never needed to calculate the dimension
of a face since the way they were produced guaranteed their dimension.
However one may wonder how one can efficiently calculate the dimension
of a face.  One of the most efficient methods makes use of the fact, discussed
in \cite{BS}, 
that the dimension is equal to $e+k-2n$, where $e$ is the number of
1's in the matrix of $F$ and $k$ the number of components in the graph
corresponding to $F$. ({\it i.e.,\/}\ the bipartite graph on $2n$
letters in which $i$ is joined to $j$ when the $(i,j)$ entry of the
matrix of $F$ is 1.) 

\item The relative volume of any $d$-face $F$ can be computed
in several different ways since it is the sum of the relative volumes
of the facets opposite any vertex of $F$.  This yields linear relations on the
volumes of $(d-1)$-faces of $B_n$.  It seems conceivable that
these linear relations could be strong enough to yield useful
information about the volumes.  However from our limited investigation
this does not appear to save anything in our computations.

\item Since our standard triangulations all involve minimum volume simplices,
one might wonder whether all minimum volume simplices with vertices
from the vertex set of $B_n$ belong to one of these triangulations.
For $n=4$, we found that there are 658584 minimum volume simplices whose
vertices are vertices of $B_4$.  Of these, only 641112 belong to some
standard triangulation.

\end{itemize}

\section{The Magic Squares Method}
\label{sec:5}

In the next two sections we describe the magic squares method for 
calculating the volume of $B_n$.  We have no reason to
believe that our implementation is substantially different 
from those used by others.  (See  \cite{DG}, \cite{Mo},\cite{SS}, 
and \cite{Stu}.)  The only apparent novelty
is that we have carried out the computation when $n=8$.

We briefly explain here the connection between magic squares and
the volume of $B_n$.

It is known that for a $d$-dimensional polytope $P$ with integer vertices, for
any nonnegative integer $t$, 
the number $e(P,t)$ of lattice points contained in $t\cdot P$ is a polynomial
of degree $d$ in $t$.  This polynomial is called the
{\it Ehrhart polynomial\/}\ of $P$.  Its leading coefficient is the 
volume of $P$ in units equal to the volume of the fundamental domain of the
affine lattice spanned by $P$.
Thus if we know the values of $e(P,t)$ for values of $t$ from 0 to $d$, 
we can find the Ehrhart polynomial by interpolation and in that way
determine the volume of $P$.

For $P=B_n$, this method is particularly attractive since the polynomial
is known to have certain symmetries, which make it necessary to calculate
the values of $e(B_n,t)$ for $t$ only up to and including
$\binom{n-1}{2}$ rather than $(n-1)^2$. 

Note that $e(B_n,t)$ is exactly the number of $n \times n$ matrices
with nonnegative integer entries and all row and column sums equal to $t$,
i.e., the number of $n \times n$ magic squares with sum $t$.
In the next section we will describe how to count magic squares
relatively efficiently.

To see that we need only find $e(B_n,t)$ for values of $t$ up to and
including $\binom{n-1}{2}$ we refer to the following identities:
\begin{enumerate}
\item $e(B_n,t)=0$ for $-n+1 \le t \le -1$.
\item $e(B_n,-n-t)=(-1)^{n-1} e(B_n,t)$ for all $t$.
\end{enumerate}
These identities (conjectured in \cite{ADG}) are easy consequences of 
Ehrhart's Law of Reciprocity, 
which states that, for a $d$-dimensional polytope $P$ with integer vertices, and
$t>0$, 
\[
e^*(P,t)=(-1)^d e(P,-t)
\]
where $e^*(P,t)$ denotes the number of integer points in
the interior of $P$.  See \cite[Chapter 9]{H}, and \cite{E} 
for proof and references. 

\noindent {\it Proof of 1\/}:\\
$e^*(B_n,t)$ is the number of $n\times n$ matrices with positive integer 
entries and all row and column sums equal to $t$.
Since all the entries are $\ge 1$, each row and column sum must be $\ge n$,
so $e^*(B_n,t)=0$ for $1 \le t \le n-1$.
By Ehrhart's Law of Reciprocity this implies $e(B_n,t)=0$ for
$-n+1 \le t \le -1$.

\noindent {\it Proof of 2\/}:\\
There is a one-to-one correspondence between $n\times n$ matrices with
nonnegative integer entries and  row and column sums $t$ and
$n\times n$ matrices with positive integer entries and row and column
sums $n+t$.  (Simply add 1 to each entry in matrices of the first type.)
Thus $e(B_n,t)=e^*(B_n,n+t)$.  Applying Ehrhart's Law of Reciprocity,
the right-hand-side equals $(-1)^{(n-1)^2}e(B_n,-n-t)$, which simplifies
to $(-1)^{n-1}e(B_n,-n-t)$.

\medskip
The effect of the first identity is that we know $n-1$ zeroes of $e(B_n,t)$.
We also have $e(B_n,0)=1$.  For each $t>0$, if we calculate the value of
$e(B_n,t)$, by the second identity we obtain also the value of $e(B_n,-n-t)$.
Thus if we calculate $e(B_n,t)$ for $t$ up to $\binom{n}{2}$, we have a total
of $n-1+1+2\binom{n}{2}=(n-1)^2+1$ values of the $e(B_n,t)$ so we have
enough data to find the polynomial $e(B_n,t)$ by interpolation.

\section{Counting Magic Squares}
\label{sec:6}
We now describe the method we used for counting the number of 
$n \times n$ magic squares of row and column sum $t$ for $t\le\binom{n-1}{2}$.
This seems no different from the methods used by others \cite{DG}
to carry out the smaller cases.

Given an $m$-tuple $r=(r_1,\dots,r_m)$ and an $n$-tuple
$c=(c_1,\dots,c_n)$ of nonnegative integers, we denote by $N(r,c)$
the number of nonnegative integer matrices with row sums
$r_1,\dots,r_m$ and column sums $c_1,\dots,c_n$.

There are a few computational principles.  
The first is that $N(r,c)=0$ unless $|r|=\sum_ir_i=\sum_jc_j=|c|$.
Next
note that
$N(r,c)$ is invariant under permutation of either the
$r$'s or the $c$'s.  Finally the principle that leads to
substantial computational savings is that, for any integer $k$ (usually
near $m/2$)
\[
N(r,c)=\sum_x N((r_1,\dots,r_k),x)N((r_{k+1},\dots,r_n),c-x)
\]
where the sum is over all nonnegative $n$-tuples $x$ such that $|x|=r_1+\cdots+r_k$
and $x_i \le c_i$, $i=1,\dots,n$.  This formula results from classifying
the matrices counted by $N(r,c)$ according to 
the column sums of the submatrix formed from the first $k$ rows.
For fixed column sums $x_1,\dots,x_n$, the column sums of the
submatrix formed by the remaining rows must be $c_i-x_i$.  
The total number of matrices in the class corresponding
to $x$ is the number of ways of choosing the top
submatrix multiplied by the number of ways of choosing the bottom.

The counting of  magic squares amounts to the calculation of
$N(r,c)$ with the ``constant'' $n$-tuples $r=c=(t,\dots,t)$.
For this special case there are a few simplifications.
We discuss the case when $n$ is even.  The same ideas apply
with slight modification when $n$ is odd.

Suppose that $n=2m$ and we wish to calculate $e(B_n,t)$.
From our general principle we have
\[
e(B_n,t)=\sum_y N(R,y)N(R,T-y)
\]
where $R$ is the $m$-tuple of all $t$'s, $T$ is
the $n$-tuple of all $t$'s, and 
$y$ runs over all nonnegative $n$-tuples satisfying $|y|=mt$,
and $y_i \le t$ for all $i$.  For a $k$-tuple $y=(y_1,\dots,y_k)$,
let us denote by $M(y)$ the number of distinct $k$-tuples which
arise by permuting the $y_i$'s.  So, if $z_1,\dots,z_l$ are distinct,
and $y$ is a $k$-tuple consisting of $k_1$ $z_1$'s, $k_2$ $z_2$'s, etc.,
then $M(y)=k!/(k_1!\dots k_l!)$.  In terms of this notation 
a more computationally efficient version of the preceding equation is
\begin{equation}\label{eq:a}
e(B_n,t)=\sum_y M(y)N(R,y)N(R,T-y)
\end{equation}
where now we further restrict $y$ to weakly increasing $n$-tuples.

We can apply this principle again to the calculation of 
$N(R,y)$ and $N(R,T-y)$ that appear in the last formula.  We find that
\begin{equation} \label{eq:b}
  N(R,y)=\sum_x M(x)N(x,(y_1,\dots,y_m))N(R-x,(y_{m+1},\dots,y_n)) 
\end{equation}
where now $x$ runs over all weakly increasing nonnegative $m$-tuples
with $|x|=y_1+\cdots+y_m$ and $x_i\le t$ for all $i$.

We can save an additional factor of 2 by noting that the 
quantities $N(R,y)$ are the same as $N(R,T-y)$ except
in a different order.  Thus if we save the former in a suitable
array, we can look up the latter ones in the array rather
than computing them.

Notice that the ingredients for calculating the sums 
$N(R,y)$ and $N(R,T-y)$ are the quantities
$N(x,y)$ where $x$ and $y$ vary over weakly increasing nonnegative $m$-tuples with
$x_i,y_i \le \binom{n-1}{2}$.
Thus it is sensible to precompute these
quantities and save the results before forming the sums for
$N(R,y)$ or the sum for $e(B_n,t)$.

For example, for $n=8$ we need to precompute the quantities $N(x,y)$
where $x$ and $y$ have length 4. Again it is easier to calculate
\[
N(x,y)=\sum N((x_1,x_2),z) N((x_2,x_3),y-z)
\]
where the sum is over all 4-long vectors $z$ with
$|z|=x_1+x_2$ and $z_i \le y_i$ for all $i$.  However we do not have available
the additional simplification to a sum over increasing sequences
$z$.
Thus on the right side we require the values $N(x,y)$, for pairs $(x,y)$
where $x$ has length 2 and $y$ has length 4, not necessarily weakly increasing,
where the components of
$x$ and $y$ vary up to 21.  It would be possible to precompute all the
needed values and save these as well for later use.  This might
be advantageous since these results are used several times each.
However, for simplicity, we use a subroutine to compute these, in effect
repeating the calculation of any $N(x,y)$ whenever needed.
This subroutine in turn calls
a subroutine for counting $2 \times 2$ matrices with prescribed row and column
sums which calculates 
$N((x_1,x_2),(y_1,y_2))=\min(x_1,x_2,y_1,y_2)+1$
whenever $|x|=|y|$.

The precalculation for $n=8$ requires about 20 minutes 
on a 500Mhz DEC alpha.  The remaining calculation 
also takes about 20 minutes.  The first part can be calculated
in single precision.  In the remaining parts we need some
sort of multiple precision method.  We perform the calculation
modulo several primes and combine the results with the
Chinese Remainder Theorem.  A similar program for $n=7$ requires
38 seconds.

Here are the Ehrhart polynomials $e(B_n,t)$, for $n=1,\dots,8$.
For each $n$, the coefficient of the last binomial coefficient
in the expression for $e(B_n,t)$ is the relative volume of $B_n$.
We express the Ehrhart polynomial of $B_n$ as an integer combination of
binomial coefficients 
$C(t+n-1+k,n-1+2k)=\binom{t+n-1+k}{n-1+2k}$, 
$k=0,\dots,\binom{n-1}{2}$,
because they are a basis for the polynomials satisfying
$p(-n-t)=(-1)^{n-1}p(t)$.

\begin{eqnarray*}
e(B_1,t) &=& C(t,0)\\
e(B_2,t) &=& C({t+1},1)\\
e(B_3,t) &=& C({t+2},2)+3C({t+3},4)\\
e(B_4,t) &=& C({t+3},3)+20C({t+4},5)+152C({t+5},7)+352C({t+6},9)\\
e(B_5,t) &=& C({t+4},4)+115C({t+5},6)+5390C({t+6},8)+\\
&& 101275C({t+7},10)+858650C({t+8},12)+\\
&& 3309025C({t+9},14)+4718075C({t+10},16)\\
e(B_6,t) &=& C({t+5},5)+714C({t+6},7)+196677C({t+7},9)+\\
&&  18941310C({t+8},11)+809451144C({t+9},13)+\\
&&  17914693608C({t+10},15)+223688514048C({t+11},17)+\\
&& 1633645276848C({t+12},19)+6907466271384C({t+13},21)+\\
&& 15642484909560C({t+14},23)+14666561365176C({t+15},25)\\
\end{eqnarray*}

\begin{eqnarray*}
e(B_7,t) &=& C({t+6},6)+ 5033C({t+7},8)+\\
&&  9090305C({t+8},10)+ 4562637436C({t+9},12)+\\
&&876755512997C({t+10},14)+ 80592643025748C({t+11},16)+\\
&&   4085047594855542C({t+12},18)+\\
&& 125166504299043921C({t+13},20)+\\
&&   2460507569635629206C({t+14},22)+\\
&&  32199612314177385616C({t+15},24)+\\
&&     285953447105799237366C({t+16},26)+\\
&&   1727929241168643056768C({t+17},28)+\\
&&   6989369809320320632154C({t+18},30)+\\
&&  18096158896344747268932C({t+19},32)+\\
&&  27093648035077238674360C({t+20},34)+\\
&& 17832560768358341943028C({t+21},36)\\
e(B_8,t) &=& C(t+7,7)+\\
&&     40312C(t+8,9)+\\
&&     544604804C(t+9,11)+\\
&&     1572522771472C(t+10,13)+\\
&&1433860489078360C(t+11,15)+\\
&&   546197610013169408C(t+12,17)+\\
&&   104573799019751624800C(t+13,19)+\\
&&   11404657872578818785152C(t+14,21)+\\
&&   773100275338739807806336C(t+15,23)+\\
&&   34668602440014649185072000C(t+16,25)+\\
&&    1075823106306592550013512704C(t+17,27)+\\
&&   23865735845675030268755397632C(t+18,29)+\\
&&     387264682746696963082402212768C(t+19,31)+\\
&&     4666750907574155613393947915904C(t+20,33)+\\
&&     42107239094874587731729608526080C(t+21,35)+\\
&&     284859465667770778104594682157824C(t+22,37)+\\
&&     1435919936068954265096148477657088C(t+23,39)+\\
&&     5307981556350553774098942855517184C(t+24,41)+\\
&&     13958946247270195588626193027208192C(t+25,43)+\\
&&     24706461764218063045041689495950080C(t+26,45)+\\
&&     26368507913706408235698183181290240C(t+27,47)+\\
&&     12816077964079346687829905128694016C(t+28,49)\\
\end{eqnarray*}

\section{Comparison}
\label{sec:7}
We now compare the two methods described above.

The main advantage of the first method seems to be that it applies
just as well to any face of $B_n$ as it does to $B_n$ itself.
To apply the algorithm to a face $F$ of $B_n$, we simply start
at the top level with the 0-1 matrix associated to $F$ and produce
lists of canonical subfaces as before.

In the second method it is not obvious how well one could do in computing
the volume of an arbitrary face $F$ of $B_n$.  This would amount to
counting the number of magic squares with prescribed zeros and
row and column sums $t$ for possibly as many as $\dim(F)$ values of $t$.
We would not have $e^*(F,t)=0$ for $1 \le t \le n-1$, nor would we have
$e(F,t)=e^*(F,n+t)$, because of the prescribed zeros in $F$. 
For certain $F$ ({\it e.g.,\/}\ those with the same number of prescribed
zeros in every row) we would have an analogous identity, and some
automatic roots, but in general we cannot guarantee any cutdown
in the number of values of $e(F,t)$ needed to determine the polynomial. 
Furthermore in the actual counting of magic squares with certain prescribed
zeros, we would not be able to exploit the symmetries used in our algorithm
above.

The second method however has the advantage that, for computing
volumes of $B_n$ itself, it is much more feasible in terms of memory.

The second method also computes the Ehrhart polynomial.
It seems possible that the first method could be modified to
compute Ehrhart polynomials of the faces as well as just their
volumes.  We would need to keep track of the numbers of 
simplices of each dimension in a standard triangulation instead
of just the simplices of the largest dimension.
\newpage

\section*{Appendix: Ratio of Volumes of \protect\boldmath$B_n$ and \protect\boldmath$A_n$}
\label{sec:appendix}

Consider the linear mapping $\cal L$ from $(n-1) \times (n-1)$ matrices to $n \times n$ 
matrices which sends matrix $E_{i,j}$ which is all zero except for a 1 at $(i,j)$ to
the matrix $F_{i,j}$ which is all zero except for 1's at $(i,j)$ and $(n,n)$ and 
$-1$'s at $(i,n)$ and $(n,j)$.  If we follow $\cal L$ by the addition of 
the $n \times n$ matrix that has the block form 
\[
\left(\begin{array}{ccc} 0 & J_{n-1,1}  \\ J_{1,n-1} & 2-n \end{array}\right)
\]
where $J_{k,l}$ is the all $k \times l$ matrix of all 1's, then we obtain
the affine mapping which sends $A_n$ to $B_n$.  Thus if we
denote the ratio we seek by
$R$, we find that $R^2$ is 
the determinant of the $(n-1)^2 \times (n-1)^2$ matrix
of dot products of $F_{i,j} \cdot F_{k,l}=2^{\delta_{ik}}2^{\delta_{jl}}$.  

But in general if $ x_{ik}$ and $y_{jl}$ are two $m \times m$ matrices
and $z$ is the $m^2 \times m^2$ tensor product matrix indexed by pairs $ij$
and $kl$ given by $ z_{ij,kl}=x_{ik}y_{jl} $ then $\det z = (\det x
\det y)^m $.

Our case is the special case that 
$x=y=J_{n-1,n-1}+I_{n-1}$.  Since the characteristic polynomial
of $-J_m$ is $\lambda^{m-1}(\lambda+m)$, the determinant of $J_{n-1}+I_{n-1}$ is
$n$.  It follows that $R=n^{n-1}$.


\end{document}